\documentclass[1p,10pt]{elsarticle}

\usepackage{amssymb}

\newtheorem{thm}{Theorem}
\newtheorem{lem}[thm]{Lemma}

\def\R{\mathbb R}

\newcommand{\fdem}{ \hfill $\square$ }

\journal{Appl. Math. Lett.}

\begin{document}

\begin{frontmatter}

\title{Rapid travelling waves in the nonlocal Fisher equation connect two unstable states}

\author[label1]{Matthieu Alfaro}
\author[label2]{J{\'e}r{\^o}me Coville}
\address[label1]{I3M, Universit\'e Montpellier 2, CC051, Place Eug\`ene Bataillon,  34095 Montpellier Cedex 5, France.}
\address[label2]{Equipe BIOSP, INRA Avignon, Domaine Saint Paul, Site Agroparc, 84914 Avignon Cedex 9, France.}

\begin{abstract}In this note, we give a positive answer to
a question addressed in \cite{Nad-Per-Tan}. Precisely we prove
that, for any kernel and any slope at the origin, there do exist
travelling wave solutions (actually those which are \lq\lq
rapid\rq\rq) of the nonlocal Fisher equation that connect the two
homogeneous steady states 0 (dynamically unstable) and 1. In
particular this allows situations where 1 is unstable in the sense
of Turing. Our proof does not involve any maximum principle
argument and applies to kernels with {\it fat tails}.
\end{abstract}

\begin{keyword}
Integro-differential equation \sep travelling waves \sep Turing
instability.

\MSC 45K05 \sep 35K57 \sep 35C07.
\end{keyword}

\end{frontmatter}


\section{Introduction}\label{s:intro}

In this work, we consider the nonlocal Fisher-KPP equation
\begin{equation}\label{edp}
\partial _t u=\partial _{xx} u+\mu u(1-\phi
* u)\qquad x\in \R\,,\quad t>0\,,
\end{equation}
where $\phi*u(x):=\displaystyle \int _\R u(x-y)\phi (y)\,dy$,
$\phi$ is a given smooth kernel such that
\begin{equation}\label{noyau}
\phi \geq 0\,,\quad \phi(0)>0\,,\quad \int _\R \phi =1\,,\quad
\int _{\R} z^2 \phi(z)\,dz<\infty\,,
\end{equation}
and $\mu >0$ is identified as the \lq\lq slope at the
origin\rq\rq. We are interested in travelling waves solutions
supported by the integro-differential equation (\ref{edp}). We are
therefore looking after a speed $c \in \R$ and a smooth and
bounded $u(x)$ such that
\begin{equation}\label{pb-tw}
-cu'=u''+\mu u(1-\phi *u)\qquad x\in\R\,,
\end{equation}
supplemented with the expected boundary conditions
\begin{equation}\label{expected}
u(-\infty)=1\,,\qquad u(+\infty)=0\,,
\end{equation}
or, when necessary, the weaker boundary conditions
\begin{equation}\label{weak}
\liminf_{x\to -\infty} u(x)>0\,,\qquad u(+\infty)=0\,.
\end{equation}

\vskip 3pt  \noindent{\bf Local Fisher-KPP equation.} If the
kernel $\phi$ is replaced by the Dirac $\delta$-function, then
(\ref{edp}) reduces to
\begin{equation}\label{edp-classical}
\partial _t u=\partial _{xx} u+\mu
u(1-u)\,,
\end{equation}
namely the classical Fisher-KPP equation \cite{Fis},
\cite{Kol-Pet-Pis}, for which $u\equiv 0$ is unstable and $u\equiv
1$ is stable. It is commonly used in the literature to model
phenomena arising in population genetics, or in biological
invasions. It is well known that the classical Fisher-KPP equation
admits monotonic travelling wave solutions with the expected
boundary conditions (\ref{expected}), for some semi-infinite
interval $[c^*:=2\sqrt \mu,\infty)$ of admissible wave speeds.
Moreover, these waves describe the long time behavior of solutions
of (\ref{edp-classical}) with compactly supported initial data or
initial data with exponential decay.

\vskip 3pt \noindent{\bf Nonlocal Fisher-KPP equation.} Let us
turn back to the integro-differential equation (\ref{edp}). In
population dynamics models, one can see the nonlinear term as the
intra-specific competition for resources. Its nonlocal form
indicates that individuals are competing with all other
individuals, whatever their positions. For more details on
nonlocal models, we refer to \cite{Gou} and the references
therein.

Again the uniform steady states of (\ref{edp}) are $u\equiv 0$ and
$u\equiv 1$. Nevertheless, because of the nonlocal effect, the
steady state 1 can be Turing unstable. In particular, this happens
when the Fourier transform $\hat \phi$ changes sign and $\mu$ is
large \cite{Gen-Vol-Aug}, \cite{Apr-Bes-Vol-Vou}. The situation is
therefore in contrast with (\ref{edp-classical}). Hence, for
travelling waves to be constructed, the authors in
\cite{Ber-Nad-Per-Ryz} have to ask not for the the expected
behavior $u(-\infty)=1$, but for the weaker condition
$\liminf_{-\infty}u>0$. More precisely, they prove the following.

\begin{lem}[Travelling waves constructed in
\cite{Ber-Nad-Per-Ryz}]\label{lem:nadin} For all $c\geq
c^*:=2\sqrt \mu$, there exists a travelling wave $(c,u)\in \R
\times C_b^2(\R)$ solution of (\ref{pb-tw}), with $u>0$ and with
the weak boundary conditions (\ref{weak}).

Moreover, there is $\mu _0>0$ such that, for all kernel $\phi$,
all $0<\mu<\mu_0$, these waves actually satisfy $u(-\infty)=1$.

Also, if the Fourier transform $\hat \phi$ is positive everywhere,
then, for all $\mu >0$, these waves actually satisfy
$u(-\infty)=1$.
\end{lem}

When $\hat \phi$ takes negative values and when $\mu>0$ is large,
such results do not precise if the waves can approach the Turing
unstable state 1 as $x\to -\infty$. By using numerical
approximation, the authors in \cite{Nad-Per-Tan} observe such
waves for the compactly supported kernel $\frac 1 2
\mathbf{1}_{[-1,1]}$. It should be noted that for kernels with
exponential decay one may use maximum principle arguments and then
derive some monotonicity properties. Hence, it is proved in
\cite{Fan-Zha} that waves which are rapid enough are monotone and
then approach 1 as $x\to -\infty$. Here, we allow kernels with
{\it fat tails} which are quite relevant in applications.
Precisely, we only assume that the second moment of $\phi$ is
finite. The result of this note is to prove that, even for such
kernels, rapid waves always connect 1 in $-\infty$. It reads as
follows.

\begin{thm}[Rapid waves connect two unstable states]\label{th}
Define
$$
\overline c=\overline c(\phi,\mu):=\mu \left(\int
_{\R}z^2\phi(z)\,dz\right)^{1/2} \left(\int _{\R}\phi (z)(1-\mu
\frac {z^2}2)_+\,dz\right)^{-1}\,.
$$
Then the waves constructed in Lemma \ref{lem:nadin} with speed $c
> \overline c$ actually satisfy $u(-\infty)=1$.
\end{thm}

\vskip 3pt Our proof does not use any maximum principle argument
neither any monotonicity property of the wave. Therefore, Theorem
\ref{th} allows the possibility of non monotonic waves. It leans
on $L^2$ estimates proved in Section \ref{s:second-moment}. In
Section \ref{s:conclusion}, we complete the proof of Theorem
\ref{th} and conclude with remarks on the bistable case.

\section{Investigating the behavior in $\pm \infty$}\label{s:second-moment}

This section contains the main contribution of the present note.
By rather elementary $L^2$ estimates, we find a sufficient
condition for a solution of (\ref{pb-tw}) to converge to 0 or 1 in
$-\infty$ and $+\infty$. In the sequel, for $i=1,2$, we define the
$i$-th moment of the kernel $\phi$ by
\begin{equation}\label{def:moments}
m_i:=\int _{\R} |z|^i\phi(z)\,dz\,.
\end{equation}

\begin{lem}[Sufficient condition for $u'\in L^2$]\label{lem:derivee-infini} Let $c\in
\R$ and $u\in C^2_b(\R)$ be a solution of (\ref{pb-tw}). Assume
$\mu \sqrt {m_2} \Vert u \Vert _{L^\infty} <|c|$. Then $u' \in
L^2(\R)$ and $\lim_{\pm\infty}u'=0$.\end{lem}

\noindent {\bf Proof.} Let us define $M:=\Vert u\Vert _{L^\infty}$
and $M':=\Vert u' \Vert _{L^\infty}$. Denote by $W$ a potential
associated with the underlying (local) monostable nonlinearity
i.e. $W'(x)=x(1-x)$. We rewrite the equation as
$$
cu'=-u''-\mu u(1-u)-\mu u(u-\phi *u)\,, $$
 multiply it by $u'$, and then
integrate from $-A<0$ to $B>0$ to get
\begin{equation}\label{eq1}
c\int _{-A}^B {u'}^2=\left[-\frac 12 {u'}^2-\mu
W(u)\right]_{-A}^B-\mu \int_{-A}^B u'u(u-\phi*u)\,.
\end{equation}
We denote by $I_{A,B}$ the last integral appearing above and use
the Cauchy-Schwarz inequality to see
\begin{equation}\label{eq2}
{I_{A,B}}^2\leq \int _{-A}^B (u'u)^2 \int _{-A}^B (u-\phi *u)^2
\leq M^2 \int _{-A}^B {u'}^2 \int _{-A}^B (u-\phi*u)^2\,.
\end{equation}
Now, for a given $x$, we write
$$
(u-\phi*u)(x)=\int _{\R} \phi(x-y)(u(x)-u(y))\,dy=\int _{\R} \int
_0 ^1 \phi (x-y) (x-y)u'(x+t(y-x))\,dtdy\,,
$$
so that another application of the Cauchy-Schwarz inequality
yields
\begin{eqnarray*}
(u-\phi*u)^2(x)& \leq &\int_{\R} \int _0^1
\phi(x-y)(x-y)^2\,dtdy\;
\int_{\R}\int_0 ^1 \phi(x-y){u'}^2(x+t(y-x))\,dtdy\\
&\leq& m_2 \int _0 ^1 \int_{\R} \phi (-z){u'}^2(x+tz)\,dzdt\,.
\end{eqnarray*}
Integrating this we discover
$$
\int_{-A}^B (u-\phi*u)^2 \leq m_2 \int _0^1\int_ {\R} \phi(-z)
\int_{-A+tz}^{B+tz}{u'}^2(y)\,dydzdt\,.
$$
Since $|u'|\leq M'$ we get, by cutting into three pieces,
$$
\int_{-A+tz}^{B+tz}{u'}^2\leq  \int_{-A}^{B}{u'}^2+2t|z|M'^2\,, $$
which in turn implies
\begin{eqnarray}
\int_{-A}^B (u-\phi*u)^2&\leq& m_2 \int_{-A}^{B}{u'}^2 +2m_2 {M'}
^2
\int _0 ^1 \int _{\R} t\phi(-z)|z|\,dzdt\nonumber\\
&\leq& m_2 \int_{-A}^{B}{u'}^2+m_2m_1M'^2\,.\label{eq4}
\end{eqnarray}
If $R_{A,B}:=\int_{-A}^{B}{u'}^2$, combining (\ref{eq1}),
(\ref{eq2}) and (\ref{eq4}) we see that
\begin{eqnarray*}
|c| R_{A,B}&\leq& \left|\left[-\frac 12 {u'}^2-\mu
W(u)\right]_{-A}^B\right|
+\mu M\sqrt {R_{A,B}}\sqrt{m_2 R_{A,B}+m_2 m_1 M'^2}\\
&\leq& M'^2+2\mu \Vert W\Vert_{L^\infty(-M,M)}+\mu \sqrt{m_2}M
\sqrt{{R_{A,B}}^2+ m_1 {M'} ^2 R_{A,B}}\,.
\end{eqnarray*}
If $\mu \sqrt {m_2} M <|c|$ then the upper estimate compels
$R_{A,B}=\int_{-A}^B {u'}^2$ to remain bounded, so that $u'\in
L^2$. Since $u'$ is uniformly continuous on $\R$, this implies
$\lim _{\pm \infty}u'=0$. This concludes the proof of the lemma.
 \fdem

\begin{lem}[Sufficient condition for $u(\pm\infty)\in \{0,1\}$]\label{lem:sufficient} Let $c\in
\R$ and $u\in C^2_b(\R)$ be a solution of (\ref{pb-tw}). Assume
\begin{equation}\label{condition}
\mu \sqrt {m_2} \Vert u \Vert _{L^\infty} <|c|\,.
\end{equation}
Then $\lim _{+\infty} u$ and $\lim _{-\infty} u$ exist and belong
to $\{0,1\}$.\end{lem}

\noindent {\bf Proof.} Since the proof is similar on both sides we
only work in $+\infty$. Denote by $\mathcal A$ the set of
accumulation points of $u$ in $+\infty$. Since $u$ is bounded,
$\mathcal A$ is not empty. Let $\theta \in \mathcal A$. There is
$x_n \to +\infty$ such that $u(x_n)\to \theta$. Then
$v_n(x):=u(x+x_n)$ solves
$$
{v_n}''+c{v_n}'=-\mu v_n(1-\phi*v_n)\quad \textrm{ on } \R\,.
$$
Since the $L^\infty$ norm of the right hand side member is
uniformly bounded with respect to $n$, the interior elliptic
estimates imply that, for all $R>0$, all $1<p<\infty$, the
sequence $(v_n)$ is bounded in $W^{2,p}([-R,R])$. From Sobolev
embedding theorem, one can extract $v_{\varphi (n)} \to v$
strongly in $C^{1,\beta}_{loc}(\R)$ and weakly in $W^{2,p}_{loc}
(\R)$. It follows from Lemma \ref{lem:derivee-infini} that
$$
v'(x)=\lim _{n\to\infty} u'(x+x_{\varphi(n)})=0\,.
$$
Combining this with the fact that $v$ solves
$$v''+cv'=-\mu v(1-\phi
*v)\quad \textrm{ on } \R\,,
$$
we collect $v\equiv 0$ or $v\equiv 1$. From $v(0)=\lim _n
u(x_{\varphi(n)})=\theta$ we deduce that $\theta\in\{0,1\}$. Since
$u$ is continuous, $\mathcal A$ is connected and therefore
$\mathcal A=\{0\}$ or $\mathcal A=\{1\}$. Therefore $u(+\infty)$
exists and is equal to $0$ or $1$.\fdem

\section{Conclusion}\label{s:conclusion}

\subsection{Proof of Theorem \ref{th}}\label{ss:proof}

Let us consider a travelling wave $(c,u)$ as in Lemma
\ref{lem:nadin}. In view of $\liminf_{-\infty}u>0$ and Lemma
\ref{lem:sufficient}, for $u(-\infty)=1$ to hold it is enough to
have (\ref{condition}). Therefore we need to investigate further
the bound $\Vert u\Vert _{L^\infty}$. To construct $(c,u)$, the
authors in \cite{Ber-Nad-Per-Ryz} first consider the problem in a
finite box $(-a,a)$. They prove {\it a priori} bounds for
solutions in the box, use a Leray-Schauder degree argument to
construct a solution $(c_a,u_a)$ in the box and then pass to the
limit $a\to \infty$ to construct $(c,u)$ a solution on the line
$\R$. One of the crucial {\it a priori} estimate is the existence
of a constant
$$ K_0=K_0(\phi,\mu):=\left(\int _{\R}\phi (z)(1-\mu \frac
{z^2}2)_+\,dz\right)^{-1} $$ such that $\Vert u\Vert
_{L^\infty}\leq K_0$ (see \cite[Lemma 3.1 and Lemma
3.10]{Ber-Nad-Per-Ryz}). Hence, any wave with speed $c
>\overline c=\mu \sqrt {m_2} K_0$ will satisfy (\ref{condition}). This completes the
proof of Theorem \ref{th}. \fdem

\subsection{Comments on the bistable case}\label{ss:bistable}

Let us conclude with a few comments concerning the bistable case.
The local equation is given by $\partial _t u=\partial _{xx} u+
u(u-\alpha)(1-u)$, where $0<\alpha<1$. It is well known that there
is a unique (up to translation) monotonic travelling wave with the
expected boundary conditions (\ref{expected}).

A difficult issue is now to search for (nontrivial) travelling
waves solutions supported by the integro-differential equation
\begin{equation}\label{edp-bistable}
\partial _t u=\partial _{xx} u +u(u-\alpha)(1-\phi *u)\qquad x\in \R\,,\quad t>0\,,
\end{equation}
that is $(c,u)$ such that
\begin{equation}\label{pb-tw-bistable}
-cu'=u''+u(u-\alpha)(1-\phi *u)\qquad x\in\R\,,
\end{equation}
supplemented with {\it ad hoc} boundary conditions. As far as we
know, no result exists for such waves. For instance, among other
things, nonlinearities such as $u(\phi*u-\alpha)(1-u)$ are treated
in \cite{Wan-Li-Ru}, but equation (\ref{edp-bistable}) does not
fall into \cite[equation
 (1.6)]{Wan-Li-Ru}. Indeed $g(u,v)=u(u-\alpha)(1-v)$ does not satisfy
 \cite[hypotheses (H1)--(H2)]{Wan-Li-Ru} (which would ensure the stability of
 $u\equiv 0$ and $\equiv 1$).
Moreover, the standard construction scheme used in
\cite{Ber-Nad-Per-Ryz} for the nonlocal Fisher-KPP equation cannot
be applied straightforwardly to the bistable case. More precisely,
proceeding as in \cite{Ber-Nad-Per-Ryz}, one can construct an
approximated solution $(c_a,u_a)$ defined on a bounded box
$(-a,a)$ and then try to pass to the limit as  $a\to \infty$.
However, the change of sign of the nonlinearity around $\alpha$
generates some difficulties in the establishment of sharp
\textit{a priori} estimates on $(c_a,u_a)$ and  makes very
delicate the comprehension of the behavior of a solution $u$ in
$\pm\infty$. In particular, the Harnack type arguments crucially
used in \cite{Ber-Nad-Per-Ryz} fail in this situation.

Nevertheless, even if the construction of travelling waves for the
bistable case is still to be addressed, the $L^2$ estimates of
Section \ref{s:second-moment} turn out to be of independent
interest: it is straightforward to derive the following analogous
of Lemma \ref{lem:sufficient} for the bistable equation
(\ref{pb-tw-bistable}).

\begin{lem}[Sufficient condition for $u(\pm\infty)\in \{0,\alpha,1\}$, bistable case] Let $c\in
\R$ and $u\in C^2_b(\R)$ be a solution of (\ref{pb-tw-bistable}).
Assume
\begin{equation}\label{condition-bistable}
 \sqrt {m_2} \Vert u \Vert _{L^\infty} ^2 <|c|\,.
\end{equation}
Then $\lim _{+\infty} u$ and $\lim _{-\infty} u$ exist and belong
to $\{0,\alpha,1\}$.\end{lem}

\vskip 3pt \noindent \textbf{Acknowledgements.} The authors would
like to thank the anonymous referees whose valuable comments led
to a real improvement of the setting of the result.

The first author is supported by the French {\it Agence Nationale
de la Recherche} within the project IDEE (ANR-2010-0112-01).


\begin{thebibliography}{ABCD}

\bibitem{Apr-Bes-Vol-Vou} N. Apreutesei, N. Bessonov, V. Volpert and V. Vougalter,
{\it Spatial structures and generalized travelling waves for an
integro-differential equation}, Discrete Contin. Dyn. Syst. Ser. B
{\bf 13} (2010), no. 3, 537--557.

\bibitem{Ber-Nad-Per-Ryz} H. Berestycki, G. Nadin, B. Perthame and
L. Ryzhik, {\it The non-local Fisher-KPP equation: travelling
waves and steady states}, Nonlinearity {\bf 22} (2009), no. 12,
2813--2844.

\bibitem{Fan-Zha} J. Fang and X.-Q. Zhao, {\it Monotone wavefronts of the nonlocal Fisher-KPP equation},
 Nonlinearity {\bf 24} (2011), 3043--3054.

\bibitem{Fis} R. A. Fisher, {\it The wave of advance of advantageous genes},
Ann. of Eugenics {\bf 7} (1937), 355--369.

\bibitem{Gen-Vol-Aug} S. Genieys, V. Volpert and P. Auger, {\it
Pattern and waves for a model in population dynamics with nonlocal
consumption of resources}, Math. Model. Nat. Phenom. {\bf 1}
(2006), no. 1, 65--82.

\bibitem{Gou} S. A. Gourley, {\it Travelling front solutions of a nonlocal Fisher
equation}, J. Math. Biol. {\bf 41} (2000), no. 3, 272--284.

\bibitem{Kol-Pet-Pis} A. N. Kolmogorov, I. G. Petrovsky and N. S. Piskunov, {\it Etude de
l'\'equation de la diffusion avec croissance de la quantit\'e de
mati\`ere et son application \`a un probl\`eme biologique},
Bulletin Universit\'e d'Etat  Moscou, Bjul. Moskowskogo Gos.
Univ., 1937, 1--26.

\bibitem{Nad-Per-Tan} G. Nadin, B. Perthame and M. Tang, {\it Can a traveling wave
connect two unstable states? The case of the nonlocal Fisher
equation}, C. R. Math. Acad. Sci. Paris {\bf 349} (2011), no.
9-10, 553--557.

\bibitem{Wan-Li-Ru} Z.-C. Wang, W.-T. Li and S. Ruan,
{\it Existence and stability of traveling wave fronts in reaction
advection diffusion equations with nonlocal delay}, J.
Differential Equations {\bf 238} (2007), 153--200.

\end{thebibliography}
\end{document}